\documentclass[twoside,a4paper,12pt,leqno]{amsart}

\usepackage{verbatim,amsxtra,amssymb,amsmath,color,psfrag,graphicx}

\definecolor{mygrey}{gray}{0.70}
\definecolor{mygreen}{rgb}{0,.75,0}
\definecolor{myred}{rgb}{1,0,0}
\definecolor{orange}{rgb}{1,.5,0}

\numberwithin{equation}{part}

\newtheorem{thm}[equation]{Theorem}

\newtheorem{cor}[equation]{Corollary}
\newtheorem{lem}[equation]{Lemma}
\newtheorem{prp}[equation]{Proposition}

\theoremstyle{definition}
\newtheorem{dfn}[equation]{Definition}

\newtheorem*{nonum1}{Theorem E (Existence of conditional probabilities)}
\newtheorem*{nonum2}{Theorem C (Classification of homomorphisms of Lebesgue spaces)}

\theoremstyle{remark}
\newtheorem{rem}[equation]{Remark}

\newcommand{\thmref}[1]{Theorem~\ref{#1}}
\newcommand{\prpref}[1]{Proposition~\ref{#1}}
\newcommand{\lemref}[1]{Lemma~\ref{#1}}
\newcommand{\corref}[1]{Corollary~\ref{#1}}
\newcommand{\dfnref}[1]{Definition~\ref{#1}}
\newcommand{\remref}[1]{Remark~\ref{#1}}

\newcommand{\secref}[1]{Section~\ref{#1}}

\newcommand\act{\circlearrowleft}

\renewcommand\b{\beta}

\newcommand\C{\mathcal C}
\newcommand\cof{\mathsf{cofinite}}

\newcommand\D{\mathcal D}
\renewcommand\d{\delta}

\newcommand\free{\mathsf{free}}

\renewcommand\H{\mathbb H}
\newcommand\HBall{\operatorname{HBall}}
\newcommand\hb{\mathsf{horB}}

\newcommand\hs{\mathsf{horS}}

\newcommand\Iso{\operatorname{Iso}}

\renewcommand\L{\Lambda}
\renewcommand\l{\lambda}
\newcommand\la{\langle}

\renewcommand\o{\omega}
\newcommand\ov{\overline}

\newcommand\pt{\partial}

\newcommand\R{\mathbb R}
\newcommand\ra{\rangle}

\newcommand\s{\sigma}

\newcommand\X{\mathcal X}

\newcommand\wt{\widetilde}

\newcommand\Z{\mathbb Z}

\begin{document}

\title{Hopf decomposition and horospheric limit sets}

\author{Vadim A. Kaimanovich}
\address{Mathematics, Jacobs University Bremen, Campus Ring 1, D-28759, Bremen, Germany}
\email{v.kaimanovich@jacobs-university.de}

\thanks{This paper was written during my stay at the Erwin Schr\"odinger Institute (ESI) in Vienna
as a senior research fellow. I would like to thank ESI for the support and
excellent working conditions.}



\date{}

\begin{abstract}
By looking at the relationship between the recurrence properties of a
countable group action with a quasi-invariant measure and the structure of its
ergodic components we establish a simple general description of the Hopf
decomposition of the action into the conservative and the dissipative parts in
terms of the Radon--Nikodym derivatives of the action. As an application we
prove that the conservative part of the boundary action of a discrete group of
isometries of a Gromov hyperbolic space with respect to any invariant
quasi-conformal stream coincides (mod 0) with the big horospheric limit set of
the group.
\end{abstract}

\maketitle

\thispagestyle{empty}

Conservativity and dissipativity are, alongside with ergodicity, the most
basic notions of the ergodic theory and go back to its mechanical and
thermodynamical origins. The famous \emph{Poincar\'e recurrence theorem}
states that any invertible transformation $T$ preserving a probability measure
$m$ on a state space $X$ is \emph{conservative} in the sense that any positive
measure subset $A\subset X$ is \emph{recurrent}, i.e., for a.e. starting point
$x\in A$ the trajectory $\{T^n x\}$ eventually returns to $A$. These
definitions obviously extend to an arbitrary measure class preserving action
$G\act(X,m)$ of a general countable group $G$ on a probability space $(X,m)$.
The opposite notions are those of \emph{dissipativity} and of a
\emph{wandering set}, i.e., such a set $A$ that all its translates $gA,\,g\in
G,$ are pairwise disjoint. An action is called \emph{dissipative} if it admits
a positive measure wandering set, and it is called \emph{completely
dissipative} if, moreover, there is a wandering set such that the union of its
translates is (mod 0) the whole action space.

Our approach to these properties is based on the observation that the notions
of conservativity and dissipativity admit a very natural interpretation in
terms of the \emph{ergodic decomposition} of the action (under the assumption
that such a decomposition exists, i.e., the action space is a \emph{Lebesgue
measure space}). Let $\C\subset X$ denote the union of all the purely
non-atomic components, and let $\D=X\setminus\C$ be the union of all the
purely atomic ergodic components. We call $\C$ and $\D$ the \emph{continual}
and \emph{discontinual} parts of the action, respectively. Further, let
$\D_\free$ be the subset of $\D$ consisting of the points with trivial
stabilizers, i.e., the union of free orbits in $\D$. The restriction of the
action to the set $\C\cup(\D\setminus\D_\free)$ is conservative, whereas the
restriction to the set $\D_\free$ is completely dissipative, thus providing
the so-called \emph{Hopf decomposition} of the action space into the
conservative and completely dissipative parts (\thmref{th:hopf}).
[Historically, such a decomposition was first established in the pioneering
paper of Eberhard~Hopf \cite{Hopf30} for one-parameter groups of measure
preserving transformations.] Although this fact is definitely known to the
specialists, it rather belongs to the ``folklore'', and the treatment of this
issue in the literature is sometimes pretty confused, so that we felt it
necessary to give a clear and concise proof.

The continual part $\C$ can be described as the set of points for which the
orbitwise sum of the Radon--Nikodym derivatives of the action is infinite
(\thmref{th:equiv}(iii)). Therefore, in the case of (mod 0) free actions this
condition completely characterizes the conservative part of the action. Once
again, the specialists in the theory of discrete equivalence relations will
hardly be surprised by this result. However, to the best of our knowledge, in
spite of its simplicity it has never been formulated explicitly. A simple
consequence of this fact is the description of the conservative part of a free
action $G\act(X,m)$ as the set of points $x\in X$ with the property that
\begin{equation*} \tag{$*$}
\textsf{there exists $t=t(x)$ such that $\{g\in G:dgm/dm(x)\ge t\}$ is
infinite}
\end{equation*}
(\thmref{th:equiv}(iv)).

The latter result completely trivializes the problem of a geometric
description of the conservative part of the boundary action of a discrete
group of isometries $G$ of a Gromov hyperbolic space $\X$ with respect to a
certain natural measure class, which is our main application.

More precisely, for any boundary point $\o\in\pt\X$ and any $x,y\in X$ let
$\b_\o(x,y) = \limsup_z [d(y,z)-d(x,z)]$, where $z\in\X$ converges to $\o$ in
the hyperbolic compactification. For CAT($-1$) spaces $\b_\o$ are the usual
\emph{Busemann cocycles}, whereas in the general case the cocycle identity is
satisfied up to a uniformly bounded additive error only, so that we have to
call them \emph{Busemann quasi-cocycles}. Then one can look for a family
$\l=\{\l_x\}$ of pairwise equivalent finite boundary measures parameterized by
points $x\in\X$ (following \cite{Kaimanovich-Lyubich05} we use the term
\emph{stream} for such families) whose mutual Radon--Nikodym derivatives are
prescribed by $\b_\o$ in the sense that
\begin{equation*} \tag{$**$}
\left| \log\frac{d\l_x}{d\l_y}(\o) - D\b_\o(x,y)\right| \le C \qquad
\forall\,x,y\in\X,\,\o\in\pt\X
\end{equation*}
for certain constants $C\ge 0, D>0$. We call such a stream
\emph{quasi-conformal of dimension~$D$}. It is \emph{invariant} if
$g\l_x=\l_{gx}$ for any $g\in G,x\in\X$.

In the CAT($-1$) case any invariant quasi-conformal stream is equivalent (with
uniformly bounded Radon--Nikodym derivatives) to an invariant conformal
\emph{stream} of the same dimension, i.e., such that the logarithms of the
Radon--Nikodym derivatives are precisely proportional to the Busemann
cocycles. Given a reference point $o\in\X$, an invariant (quasi-)conformal
stream is uniquely determined just by a \emph{(quasi-)conformal} measure
$\l_o$ with the property that $\log dg\l_o/d\l_o(\o)$ is proportional to
$\b_\o(go,o)$ (up to a uniformly bounded additive error).

Coornaert \cite{Coornaert93} proved (applying the construction used by
Patterson \cite{Patterson76} in the case of Fuchsian groups) that for any
discrete group of isometries $G$ of a Gromov hyperbolic space $\X$ with a
finite critical exponent $D$ there exists an invariant quasi-conformal stream
of dimension $D$ supported by the limit set of $G$.

It is with respect to the measure class of an invariant quasi-conformal stream
$\l=\{\l_x\}$ that we study the Hopf decomposition of the boundary action. The
geometric description ($**$) of the Radon--Nikodym derivatives in combination
with criterion ($*$) immediately implies that the conservative part of the
action coincides (mod 0) with the \emph{big horospheric limit set} $\L^\hb_G$
of the group $G$, i.e., the set of points $\o\in\pt\X$ for which there exists
$t=t(\o)$ such that $\{g\in G:\b_\o(go,o)\ge t\}$ is infinite for a certain
fixed reference point $o\in X$, or, in other words, the set of points
$\o\in\pt\X$ such that a \emph{certain} horoball centered at $\o$ contains
infinitely many points from the orbit $Go$ (\thmref{th:main}).

This characterization of the conservative part of the boundary action was
first established by Pommerenke \cite{Pommerenke76} for Fuchsian groups with
respect to the visual stream on the boundary circle (although in a somewhat
different terminology). Pommerenke's argument uses analytic properties of the
Blaschke products and does not immediately carry over to the higher
dimensional situation. Sullivan \cite{Sullivan81} used a more direct
geometrical approach and proved this characterization for Kleinian groups,
again with respect to the visual stream. Actually he considered the
\emph{small horospheric limit set} $\L^\hs_G$ (also called just
\emph{horospheric limit set}; it is defined by requiring that the intersection
of \emph{any} horoball centered at $\o\in\L^\hs_G$ with the orbit $Go$ be
infinite) essentially showing that $\L^\hb\setminus\L^\hs$ is a null set. By
extending Sullivan's approach (with some technical complications) Tukia
\cite{Tukia97} proved \thmref{th:main} for Kleinian groups with respect to an
arbitrary invariant conformal stream.

Our completely elementary approach subsumes all these particular cases and
identifies the conservative part of the boundary action with the big
horospheric limit set in full generality, for an arbitrary invariant
quasi-conformal stream on a general Gromov hyperbolic space.

It is clear from looking at criterion ($*$) that the right object in the
context of studying conservativity of boundary actions is the big horospheric
limit set $\L^\hb$ rather than the small one $\L^\hs$. Nonetheless it is
plausible that $\L^\hb\setminus\L^\hs$ is a null set with respect to any
invariant quasi-conformal stream on an arbitrary Gromov hyperbolic space. This
was proved by Sullivan \cite{Sullivan81} for Kleinian groups with respect to
the visual stream, and for subgroups of a free group (again with respect to
the uniform stream) it was done in \cite{Grigorchuk-Kaimanovich-Nagnibeda07}.

We refer the reader to \cite{Grigorchuk-Kaimanovich-Nagnibeda07} for a recent
detailed study of the interrelations between various kinds of limit sets in
the simplest model case of the action of a free subgroup on the boundary of
the ambient finitely generated free group. Actually, it was this collaboration
that brought me to the issues discussed in the present article, and I would
like to thank my collaborators Rostislav Grigorchuk and Tatiana Nagnibeda for
the gentle insistence with which they encouraged my work.

\section{Structure of the ergodic components and recurrence properties}

\subsection{Lebesgue spaces} \label{sec:L}

We begin by recalling the basic properties of the \emph{Lebesgue measure
spaces} introduced by Rokhlin, see \cite{Rohlin52,Cornfeld-Fomin-Sinai82}.
Measure-theoretically these are the probability spaces such that their
non-atomic part is isomorphic to an interval with the Lebesgue measure on it.
Thus, any Lebesgue space is uniquely characterized by its \emph{signature}
$\s=(\s_0;\s_1,\s_2,\dots)$, where $\s_0$ is the total mass of the non-atomic
part, and $\s_1\ge\s_2\ge\dots$ is the ordered sequence of the values of its
atoms (extended by zeroes if the set of atoms is finite or empty). There is
also an intrinsic definition of Lebesgue spaces based on their separability
properties. However, for applications it is usually enough to know that any
\emph{Polish topological space} (i.e., separable, metrizable, complete)
endowed with a Borel probability measure is a Lebesgue measure space. We shall
follow the standard measure theoretical convention:

\medskip

\emph{Unless otherwise specified, all the identities, properties etc. related
to measure spaces will be understood {\rm mod~0} (i.e., up to null sets)}. In
particular, all the $\s$-algebras are assumed to be \emph{complete}, i.e., to
contain all the measure 0 sets.

\medskip

An important feature of the Lebesgue measure spaces is

\begin{nonum1}
Let $p:(X,m)\to (\ov X,\ov m)$ be a \emph{homomorphism} (\emph{projection,
factorization, quotient map\/}) of Lebesgue spaces, i.e., for any measurable
set $\ov A\subset\ov X$ its preimage $A=p^{-1}(\ov A)\subset X$ is also
measurable, and $m(A)=\ov m(\ov A)$. Then the preimages $X_{\ov x}=p^{-1}(\ov
x),\, \ov x\in\ov X$, can be uniquely endowed with \emph{conditional}
probability measures $m_{\ov x}$ in such a way that $(X_{\ov x}, m_{\ov x})$
are Lebesgue spaces and the measure $m$ decomposes into an integral of the
measures $m_{\ov x},\,\ov x\in\ov X,$ with respect to the \emph{quotient
measure} $\ov m$ on~$\ov X$. Namely, for any function $f\in L^1(X,m)$ its
restrictions $f_{\ov x}$ to $X_{\ov x}$ are measurable and belong to the
respective spaces $L^1(X_{\ov x},m_{\ov x})$, the integrals $\ov f(\ov x)= \la
f_{\ov x},m_{\ov x}\ra$ depend on~$\ov x$ measurably, and $\la f,m\ra = \la
\ov f, \ov m\ra$ (cf. the classical Fubini theorem).
\end{nonum1}

In fact, the above property follows from the classical Fubini theorem in view
of Rokhlin's

\begin{nonum2}
Any homomorphism $p:(X,m)\to (\ov X,\ov m)$ of Lebesgue spaces is uniquely (up
to an isomorphism) determined by the signatures of the quotient measure $\ov
m$ and of the conditional measures $m_{\ov x}$. Namely, let us denote by $I$
the unit interval endowed with the Lebesgue measure $\l$, and partition $I$
into a union of consecutive intervals $I_0,I_1,\dots$ with $\l(I_i)=\s_i$ for
a certain signature $\s$. Further, let us consider a coordinatewise measurable
assignment of signatures $\s^x$ to points $x\in I$ which is constant on the
intervals $I_1,I_2,\dots$, and, as before, let $I_0^x,I_1^x,\dots$ be the
consecutive subintervals of $I$ with $\l(I_i^x)=\s_i^x$. Denote by $(X,m)$ the
Lebesgue space obtained from the square $(I\times I,\l\otimes\l)$ by
collapsing the sets $\{x\}\times I_i^x,\,i\ge 1,x\in I_0,$ and $I_j\times
I_i^x,\,i,j\ge 1, x\in I_j,$ onto single points. In the same way, let $(\ov
X,\ov m)$ be the quotient space of the interval $(I,\l)$ obtained by
collapsing the intervals $I_1,I_2,\dots$ onto single points, so that the
signature of $(\ov X,\ov m)$ is $\s$. Then the projection of the square
$I\times I$ onto the first coordinate gives rise to a homomorphism from
$(X,m)$ to $(\ov X,\ov m)$, and the signatures of the associated conditional
measures are precisely $\s^x$. The claim is that any homomorphism of Lebesgue
spaces can be obtained in this way. In particular, if both the quotient
measure $\ov m$ and all the conditional measures $m_{\ov x}$ are purely
non-atomic (i.e., have the signature $(1;0,0,\dots)$), then the corresponding
quotient map is isomorphic just to the projection of the unit square onto the
first coordinate.
\end{nonum2}

Obviously, any homomorphism of Lebesgue spaces gives rise to the
\emph{preimage sub-$\s$-algebra} in $X$ which consists of the preimages of all
the measurable sets in $\ov X$. Another important feature of the Lebesgue
spaces is that, in fact, an arbitrary sub-$\s$-algebra in $X$ can be obtained
in this way for a certain uniquely defined quotient map.

\medskip

Below we shall use the following elementary fact which follows at once from
the uniqueness of the system of conditional measures.

\begin{lem} \label{lem:RNsame}
Let $T$ be an invertible measure class preserving transformation of a Lebesgue
space $(X,m)$, and let $p:(X,m)\to (\ov X,\ov m)$ be its $T$-invariant
projection, i.e., $p(Tx)=p(x)$ for a.e. $x\in X$. Then the conditional
measures $m_{\ov x}, \,\ov x\in\ov X,$ of the projection $p$ are
quasi-invariant with respect to $T$ and have the same Radon--Nikodym
derivatives as the measure~$m$:
$$
\frac{d\,Tm}{dm}(x) = \frac{d\,Tm_{\ov x}}{dm_{\ov x}}(x) \;, \qquad
\text{where} \quad \ov x = p(x) \;,
$$
for a.e. $x\in X$.
\end{lem}

\subsection{Ergodic components, continuality and discontinuality} \label{sec:b}

\emph{Let now $G\act(X,m)$ be an action of an infinite countable group $G$ by
measure class preserving transformations on a Lebesgue space $(X,m)$ --- which
will be our standing assumption through the rest of this Section.}

The quotient space $(\ov X,\ov m)$ of $(X,m)$ determined by the $\s$-algebra
of $G$-invariant sets is called the \emph{space of ergodic components} of the
action of $G$ on the space $(X,m)$, and the preimages $X_{\ov x}$ endowed with
the conditional measures $m_{\ov x}$ are called the \emph{ergodic components}.
The ergodic components $X_{\ov x}$ are $G$-invariant, the conditional measures
$m_{\ov x}$ are $G$-quasi-invariant, and the action of $G$ on the spaces
$(X_{\ov x}, m_{\ov x})$ is ergodic (e.g., see \cite{Schmidt77}).
\lemref{lem:RNsame} implies that the conditional measures $m_{\ov x}$ have the
same Radon--Nikodym derivatives with respect to the action of $G$ as the
original measure $m$.

Since the ergodic components are ergodic, each of them is either purely atomic
(in which case it consists of a single $G$-orbit), or purely non-atomic.

\begin{dfn} \label{def:hopf}
The \emph{continual} $\C$ (resp., \emph{discontinual} $\D$) \emph{part} of the
action $G\act(X,m)$ is the union of all the purely non-atomic (resp., purely
atomic) components of the action. Denote by $\ov\C$ and $\ov\D$ the
corresponding $G$-invariant subsets of the space of ergodic components $\ov X$
(their measurability follows from Theorem C). An orbit $Gx$ is
\emph{continual} (resp., \emph{discontinual}) if it belongs to $\C$ (resp., to
$\D$). The action is \emph{discontinual} if $m(\D)>0$ and \emph{continual}
otherwise. If $m(\C)=0$ the action is called \emph{completely discontinual}.
\end{dfn}

\begin{rem}
The quotient measure class on the space of ergodic components and the measure
classes of the conditional measures on the ergodic components do not change
when the measure $m$ is replaced with an equivalent one, so that the above
Definition (as well as various definitions below related to conservativity and
dissipativity) depends only on the measure class of $m$.
\end{rem}

\begin{lem} \label{lem:choice}
Let $A\subset X$ be a measurable $G$-invariant subset. It is contained in the
discontinual part $\D$ if and only if one can select, in a measurable way, a
representative from each $G$-orbit contained in $A$, i.e., if and only if
there exists a measurable map $\pi:A\to A$ which is constant along the orbits
of the action.
\end{lem}

\begin{proof}
If $\pi$ is such a map, then it identifies the space of ergodic components of
$A$ with $\pi(A)$, so that in particular $A\subset\D$. Conversely, by Theorem
C the discontinual part $\D$ can be identified, in a measurable way, with the
product $\ov\D\times\{1,2,\dots\}$. Then one can define the map $\pi:\D\to\D$
with the required properties as $\pi(\ov x,n)=(\ov x,1)$.
\end{proof}

\begin{rem} \label{rem:choice}
Below we shall also encounter the situation when instead of a map with the
properties from the above Lemma one has an orbit constant measurable map
$x\mapsto M_x$, where $M_x$ is a non-empty finite \emph{subset} of the orbit
$Gx$. This situation can be easily reduced to \lemref{lem:choice} by choosing
(in a measurable way!) just a single point from each of the sets $M_x$. This
can be done, for instance, by identifying the space $(X,m)$ with the unit
interval (with possible collapsing corresponding to the atoms of the measure
$m$) and taking then the minimal of the points of $M_x$.
\end{rem}

\subsection{Recurrent and wandering sets} \label{sec:rw}

Let us first remind the definitions (e.g., see \cite{Krengel85} for the case
when $G$ is the group of integers $\Z$). Actually, we have to slightly modify
them (and to distinguish recurrence from infinite recurrence) in order to take
into account certain effects which do not arise for the group $\Z$ (see below
the remarks after \thmref{th:hopf}).

\begin{dfn} \label{dfn:rec}
A measurable set $A\subset X$ is called \emph{recurrent} (resp.,
\emph{infinitely recurrent}) if for a.e. point $x\in A$ the trajectory $Gx$
eventually returns to $A$, i.e., $gx\in A$ for a certain element $g\in G$
other than the group identity $e$ (resp., returns to $A$ infinitely often,
i.e., $gx\in A$ for infinitely many elements $g\in G$). The opposite notion is
that of a \emph{wandering set} (kind of a ``fundamental domain''), i.e., a
measurable set $A\subset X$ with pairwise disjoint translates $gA,\,g\in G$.
\end{dfn}

We shall now explain the connection between these notions and our
\dfnref{def:hopf}.

\begin{prp} \label{prp:rec}
Any measurable subset of the continual part $\C$ is infinitely recurrent.
\end{prp}

\begin{proof}
For a measurable subset $A\subset\C$ put
$$
A_0 = \{x\in A: gx\in A \;\text{for finitely many}\; g\in G \} \subset F \;,
$$
where
$$
F = \{x\in X: Gx\cap A \; \text{is non-empty and finite} \} \;.
$$
The set $F$ is $G$-invariant and measurable, and by \lemref{lem:choice} and
\remref{rem:choice} $F\subset\D$. Therefore $A_0$ must be a null set, whence
the claim.
\end{proof}

Denote by $\D_\free$ (resp., $\D_\cof$) the union of all the \emph{free}
(resp., \emph{cofinite}) discontinual orbits, i.e., such that the stabilizers
of their points are trivial (resp., finite). Obviously,
$\D_\free\subset\D_\cof$, and both these sets are measurable. Let
$\ov\D_\free$ and $\ov\D_\cof$ be the corresponding subsets of the space of
ergodic components $\ov X$. The definition immediately implies

\begin{prp} \label{pr:disc}
Any measurable subset of $\D\setminus\D_\free$ (resp., of
$\D\setminus\D_\cof$) is recurrent (resp., infinitely recurrent).
\end{prp}

Let us now look at the wandering sets.

\begin{prp} \label{pr:wand}
Any wandering set is contained in $\D_\free$, and there is a maximal wandering
set, i.e., such that $\D_\free=\bigcup_{g\in G} gA$.
\end{prp}

\begin{proof}
If $A$ is a wandering set, then the map $\pi$ from the $G$-invariant union
$\wt A=\bigcup_{g\in G} gA$ to $A$ defined as $\pi(x)=Gx\cap A$ is measurable
and $G$-invariant, so that $A\subset\D$ by \lemref{lem:choice}. Moreover, all
the orbits intersecting $A$ are obviously free, whence $A\subset\D_\free$.

As for the maximality, one can take for such a wandering set any measurable
section of the projection $\D_\free\to\ov\D_\free$ (which exists by
\lemref{lem:choice}).
\end{proof}

\subsection{Hopf decomposition}

\begin{dfn}
An action $G\act(X,m)$ is called \emph{conservative} (resp., \emph{infinitely
conservative}) if any measurable subset $A\subset X$ is recurrent (resp.,
infinitely recurrent). It is called \emph{dissipative} if there is a
non-trivial wandering set, and \emph{completely dissipative} if the whole
action space $X$ is the union of translates of a certain wandering set.
\end{dfn}

The action on the disjoint union of two $G$-invariant sets is conservative
(resp., infinitely conservative) if and only if the action on each of these
sets has the same property. Taking stock of the Propositions from
\secref{sec:rw} we now obtain

\begin{thm}[Hopf decomposition for general actions] \label{th:hopf}
The action space $X$ can be decomposed into the disjoint union of two
$G$-invariant measurable sets (called its \emph{conservative} and
\emph{dissipative parts}, respectively)
$$
X = [\C \cup (\D\setminus\D_\free) ] \sqcup \D_\free
$$
such that the restriction of the action to $C\cup(\D\setminus\D_\free)$ is
conservative and the restriction to $\D_\free$ is totally dissipative.
\end{thm}

\begin{cor} \label{cor:DD}
If the action $G\act(X,m)$ is free, i.e., $\D=\D_\free$, then its conservative
part coincides with the continual set $\C$, and its dissipative part coincides
with the discontinual set $\D$.
\end{cor}

\begin{cor}[Poincar\'e recurrence theorem] \label{cor:Poin}
If the measure $m$ is invariant, then $\D_\free$ is a null set, and therefore
the action is conservative.
\end{cor}

\begin{rem}
The decomposition into the conservative and totally dissipative parts
described in the above Theorem is unique. Indeed, let $X=C_1\sqcup
D_1=C_2\sqcup D_2$ be two such decompositions. If they are different, then one
of the sets $C_1\cap D_2, C_2\cap D_1$ must be non-empty. Let it be, for
instance, $A=C_1\cap D_2$. Then the restriction of the action to $A$ has to be
simultaneously conservative (because $A\subset C_1$) and totally dissipative
(because $A\subset D_2$), which is impossible.
\end{rem}

\begin{rem}
If the conservativity is replaced with the infinite conservativity, then,
generally speaking, the Hopf decomposition as above is not possible. Namely,
the action is infinitely conservative on $\C\cup(\D\setminus\D_\cof)$ and is
completely dissipative on $\D_\free$, whereas on the remaining set
$\D_\cof\setminus\D_\free$ it is neither infinitely conservative nor
dissipative.
\end{rem}

\begin{rem}
The group $\Z$ does not contain non-trivial finite subgroups, so that for its
actions $\D_\cof=\D_\free$, and therefore infinite conservativity is
equivalent to plain conservativity. The conservative part of the action in
this case is the union of the continual part~$\C$ and the set of all the
periodic points $\D\setminus\D_\cof$.
\end{rem}

\begin{rem}
When dealing with the actions of the group $\Z$ one sometimes defines the
notion of recurrence by looking only at the ``positive semi-orbits'' $\Z_+x$.
All measurable subsets of $\C$ are recurrent in this sense as well. Indeed,
for a subset $A\subset\C$ let $A_0=\{x\in A: \Z_+x\cap A \;\text{is
finite}\}$. Then $x\mapsto zx$, where $z$ is the maximal element of $\Z$ with
$zx\in A$, is a measurable ``selection map'' in the sense of
\lemref{lem:choice}, so that $A_0\subset\D$, whence $A_0$ is a null set.
\end{rem}

\begin{rem}
If one defines \emph{strict recurrence} and \emph{strict conservativity} by
requiring that the orbit $Gx$ returns to the set $A$ at an orbit point
\emph{different} from the starting point $x$, then the strictly conservative
part of the action coincides just with the continual part $\C$. Again, as with
the infinite conservativity, the action on the set $\D\setminus\D_\free$ will
be neither strictly conservative nor dissipative.
\end{rem}

\subsection{A continuality criterion}

\begin{thm} \label{th:equiv}
Let $G\act(X,m)$ be a free measure class preserving action of a countable
group $G$ on a Lebesgue space. Denote by $\mu_x,\,x\in X$, the measure on the
orbit $Gx$ defined as
$$
\mu_x(gx)=\frac{dg^{-1} m}{dm}(x)=\frac{dm(gx)}{dm(x)}
$$
(obviously, the measures $\mu_x$ corresponding to different points $x$ from
the same $G$-orbit are proportional). Then for a.e. point $x\in X$ the
following conditions are equivalent:
\begin{itemize}
    \item[(i)]
    The orbit $Gx$ is dissipative;
    \item[(ii)]
    The orbit $Gx$ is discontinual;
    \item[(iii)]
    The measure $\mu_x$ is finite;
    \item[(iv)]
    For any $t>0$ the set $\{y\in Gx: \mu_x(y)\ge t\}$ is finite;
    \item[(v)]
    The set $M_x$ of maximal weight atoms of the measure $\mu_x$ is non-empty and
    finite.
\end{itemize}
\end{thm}

\begin{proof}
(i)$\;\Longleftrightarrow\;$(ii). This is \corref{cor:DD}.

(ii)$\implies$(iii). By \dfnref{def:hopf}, the orbit $Gx$ is discontinual if
and only if it is an ergodic component of the $G$-action on $X$. By
\lemref{lem:RNsame} in this case the measure $\pi_x$ is proportional to the
conditional measure on this ergodic component, and therefore it is finite.

(iii)$\implies$(iv)$\implies$(v). Obvious.

(v)$\implies$(i). Follows from \lemref{lem:choice} and \remref{rem:choice}
(because the map $x\mapsto M_x$ is measurable in view of Theorem C).
\end{proof}

\begin{cor} \label{cr:equiv}
Under conditions of \thmref{th:equiv} the continual part $\C$ and the
discontinual part $\D$ of the action coincide \textup{(mod~0)} with the sets
\begin{equation} \label{eq:cons}
\left\{x \in X: \sum_{g\in G}\frac{dgm(x)}{dm(x)} = \infty \right\}
\end{equation}
and
\begin{equation} \label{eq:diss}
\left\{x \in X: \sum_{g\in G}\frac{dgm(x)}{dm(x)} < \infty \right\} \;,
\end{equation}
respectively.
\end{cor}

\begin{rem}
\thmref{th:equiv} (except for the equivalence
(i)$\;\Longleftrightarrow\;$(ii)) and \corref{cr:equiv} (with the summation
taken over the equivalence class of $x$) are also true for an arbitrary
\emph{countable non-singular equivalence relation} on a Lebesgue space
$(X,m)$. In this case $\mu_x$ are the measures on the equivalence classes
determined by the \emph{Radon--Nikodym cocycle} of the equivalence relation.
In particular, the equivalence of conditions (ii), (iii), (iv), (v) from
\thmref{th:equiv} also holds for non-free actions. In what concerns
\corref{cr:equiv}, the only difference with the free case is that one has to
replace the summation over $g$ in formulas \eqref{eq:cons} and \eqref{eq:diss}
with the summation over the orbit of~$x$. The proof is precisely the same;
however, in order to spare the reader the trouble of going through the
definitions from the ergodic theory of equivalence relations (see
\cite{Feldman-Moore77}) we confine ourselves just to free actions. This
generality is sufficient, on one hand, to expose our (very simple) line of
argument, and, on the other hand, to deal with our main application to
boundary actions (\thmref{th:main}).
\end{rem}

\begin{rem} \label{rem:special}
Condition (iv) from \thmref{th:equiv} is \emph{not}, generally speaking,
equivalent just to \emph{existence} of $t>0$ such that the set $\{y\in Gx:
\mu_x(y)\ge t\}$ is finite (i.e., to boundedness of the values of the weights
of the measure $\mu_x$). The most manifest example of this is an action with a
finite invariant measure, see \corref{cor:Poin}.
\end{rem}

\section{Application to boundary actions}

\subsection{Hyperbolic spaces and limit sets}

Recall that a non-compact complete proper metric space $\X$ is \emph{Gromov
$\d$-hyperbolic} (with $\d\ge 0$) if its metric $d$ satisfies the
\emph{$\d$-ultrametric inequality}
$$
(x|z)_o \ge \min\{(x|y)_o,(y|z)_o\} - \d \qquad \forall\,x,y,z,o\in\X\;,
$$
where
$$
(x|y)_o = \frac12 \left[d(o,x) + d(o,y) - d(x,y) \right]
$$
is the \emph{Gromov product}. In addition we require that the space $X$ be
\emph{separable}. On the other hand, we do \emph{not} require the space $\X$
to be \emph{geodesic}. This class of spaces contains Cartan--Hadamard
manifolds with pinched sectional curvatures (in particular, the classical
hyperbolic spaces of constant negative curvature) and metric trees, see
\cite{Gromov87,Ghys-delaHarpe90} for more details.

A Gromov hyperbolic space $\X$ admits a natural \emph{hyperbolic
compactification} $\ov\X=\X\cup\pt\X$, and the action of the isometry group
$\Iso(\X)$ extends by continuity to a continuous \emph{boundary action} on
$\pt\X$. The boundary $\pt\X$ is a Polish space.

The \emph{limit set} $\L=\L_G\subset\pt\X$ of a discrete subgroup
$G\subset\Iso(\X)$ (any such subgroup is at most countable) is the set of all
the limit points of any given orbit $Go,\,o\in\X$, with respect to the
hyperbolic compactification, so that the closure of the orbit $Go$ in the
hyperbolic compactification is $Go\cup\L_G$ (this definition does not depend
on the choice of the basepoint $o$). The limit set is closed and
$G$-invariant. Moreover, the action of $G$ on $\L_G$ is \emph{minimal} (there
are no proper $G$-invariant closed subsets), whereas the action of $G$ on the
complement $\pt\X\setminus\L_G$ is \emph{properly discontinuous} (no orbit has
accumulation points) \cite{Gromov87,Bourdon95}.

The latter result provides a topological decomposition of the boundary action.
On the other hand, the situation is more complicated from the
measure-theoretical point of view. Let $m$ be a $G$-quasi-invariant measure on
$\pt\X$. Without loss of generality we can assume that it is purely
non-atomic. Then, since any element of $G$ fixes at most two boundary points,
the action is free. The complement $\pt\X\setminus\L_G$ is obviously contained
in the dissipative ($\equiv$ discontinual) part of the action. However, this
is as much as can \emph{a priori} be said about the ergodic properties of the
boundary action. In particular, the action on $\L_G$ need not be ergodic or
conservative. There are numerous examples witnessing to this; see
\cite{Grigorchuk-Kaimanovich-Nagnibeda07} for a detailed discussion of the
simplest model case of the action of a subgroup of a free group on the
boundary of the ambient group and for further references.

One can specialize the type of convergence in the definition of the limit set.
For instance, the \emph{radial limit set} $\L^{\textsf{rad}}$ is the set of
all the accumulation points of any fixed orbit $Go,\,o\in\X$, which stay
inside a \emph{tubular neighbourhood} of a certain geodesic ray in~$\X$. Yet
another type of the boundary convergence, which we are going to describe
below, is provided by \emph{horospheric neighborhoods}.

\medskip

Denote by
$$
 \b_z(x,y)
 = d(y,z) - d(x,z) \qquad x,y\in\X
$$
the \emph{distance cocycle} associated with a point $z\in\X$, and, following
\cite{Kaimanovich04}, put
\begin{equation} \label{eq:beta}
 \b_\o(x,y) = \limsup_{z\to\o} \b_z(x,y)
 \qquad \forall\,x,y\in\X,\,\o\in\pt\X \;.
\end{equation}
If the space $\X$ is CAT($-1$) (e.g., a Cartan--Hadamard manifold with pinched
sectional curvatures or a tree), then $\limsup$ in the above formula can be
replaced just with the ordinary limit, and $\b_\o$ are the boundary
\emph{Busemann cocycles}. Although for a general Gromov hyperbolic space
$\b_\o$ are not, generally speaking, cocycles, they still satisfy the cocycle
identity with a uniformly bounded error (i.e., they are
\emph{quasi-cocycles}). Namely,

\begin{prp}[\cite{Kaimanovich04}] \label{pr:quasicocycle}
There exists a constant $C\ge 0$ depending on the hyperbolicity constant $\d$
of the space $\X$ only such that for any $\o\in\pt\X$ the function $\b_\o$
\eqref{eq:beta} has the following properties:
\begin{itemize}
\item[(i)]
$\b_\o$ is ``jointly Lipschitz'', i.e., $|\b_\o(x,y)|\le d(x,y)$ for all
$x,y\in\X$, in particular, $\b_\o(x,x)\equiv 0\;;$
\item[(ii)]
$0 \le  \b_\o(x,y) + \b_\o(y,z) + \b_\o(z,x) \le C$ for all $x,y,z\in\X\;.$
\end{itemize}
\end{prp}

The quasi-cocycles $\b_\o$ are obviously invariant with respect to the
isometries of $\X$, i.e.,
$$
\b_{g\o}(gx,gy) = \b_\o(x,y)
\qquad\forall\,x,y\in\X,\,\o\in\pt\X,\,g\in\Iso(\X) \;.
$$

We shall define the \emph{horoball} in $\X$ centered at a boundary point
$\o\in\pt\X$ and passing through a point $o\in\X$ as
$$
\HBall_\o(o) = \{x\in\X: \b_\o(o,x)\le 0 \} \;.
$$

\begin{dfn}
The \emph{big} (resp., \emph{small}) \emph{horospheric limit set}
$\L^\hb=\L^\hb_G$ (resp., $\L^\hs=\L^\hs_G$) of a discrete group $G$ of
isometries of a Gromov hyperbolic space $\X$ is the set of all the points
$\o\in\pt\X$ such that a certain (resp., any) horoball centered at $\o$
contains infinitely many points from a fixed orbit $Go,\,o\in\X$ (the
resulting set does not depend on the choice of the orbit $Go$, see
\remref{rem:ind} below).
\end{dfn}

\begin{rem} \label{rem:ind}
As it follows from \prpref{pr:quasicocycle}, for any fixed reference point
$o\in\X$ the big (resp., small) horospheric limit sets can also be defined as
the set of all the points $\o\in\pt\X$ for which the set
$$
\{x\in Go: \b_\o(o,x) \le t \}
$$
is infinite for a certain (resp., for any) $t\in\R$.
\end{rem}

\begin{rem}
Usually our small horospheric limit set is called just the
\emph{ho\-ro\-sphe\-ric limit set}, and in the context of Fuchsian groups its
definition, along with the definition of the radial limit set, goes back to
Hedlund \cite{Hedlund36}. Following \cite{Matsuzaki02} (in the Kleinian case)
we call it \emph{small} in order to better distinguish it from the \emph{big}
one, which, although apparently first explicitly introduced by Tukia
\cite{Tukia97} (again just in the Kleinian case), essentially appears (for
Fuchsian groups) already in Pommerenke's paper \cite{Pommerenke76}.
\end{rem}

The horospheric limit sets $\L^\hs,\L^\hb$ are obviously $G$-invariant, Borel,
and contained in the full limit set $\L$ (because the only boundary
accumulation point of any horoball is just its center).

\subsection{Boundary conformal streams}

\begin{dfn}
A family of pairwise equivalent finite measures $\l=\{\l_x\}$ on the boundary
$\pt\X$ of a Gromov hyperbolic space $\X$ parameterized by points $x\in\X$ is
called a \emph{quasi-conformal stream} of dimension $D>0$ if there exists a
constant $C>0$ such that
$$
\left| \log\frac{d\l_x}{d\l_y}(\o) - D\b_\o(x,y)\right| \le C \qquad
\forall\,x,y\in\X,\,\o\in\pt\X \;.
$$
A stream $\l$ is \emph{invariant} with respect to a group $G\subset\Iso(\X)$
if
$$
\l_{gx} = g\l_x \qquad\forall\,g\in G,\, x\in\X \;.
$$
\end{dfn}

\begin{rem}
We follow here the terminology developed in \cite{Kaimanovich-Lyubich05}. More
traditionally, any invariant quasi-conformal stream is determined just by a
single finite boundary \emph{quasi-conformal measure} $\l=\l_o$ with the
property that
\begin{equation} \label{eq:cac}
\left| \frac{dg\l}{d\l}(\o) - D\b_\o(go,o) \right| \le C \qquad\forall\,g\in
G,\,\o\in\pt\X \;.
\end{equation}
for a certain reference point $o\in\X$. If $\b_\o$ are cocycles (which is the
case for CAT($-1$) spaces), then any measure $\l$ satisfying \eqref{eq:cac} is
equivalent to a unique finite measure $\l'$ (called \emph{conformal}) which
satisfies formula \eqref{eq:cac} with $C=0$ (it follows from the fact that any
uniformly bounded cocycle is cohomologically trivial). This definition is
motivated by the fact that the \emph{visual measure} on the boundary sphere
$\pt\H^{d+1}$ of the classical $(d+1)$-dimensional hyperbolic space with
sectional curvature $-1$ is conformal of dimension $d$ (in our terminology it
means that the \emph{visual stream} which consists in assigning to any point
from the hyperbolic space the associated visual measure is conformal).
However, the limit set of the group can be ``much smaller'' than the boundary
sphere and be a null set with respect to the visual measure. Existence of
conformal measures which are concentrated on the limit set (and for which the
dimension coincides with the critical exponent of the group) was first
established by Patterson \cite{Patterson76} in the case of Fuchsian groups.
His construction was further generalized (see \cite{Sullivan79},
\cite{Kaimanovich90}), ultimately providing existence of a conformal measure
for any closed subgroup of isometries of a general CAT($-1$)
space~\cite{Burger-Mozes96}. For discrete isometry groups of general Gromov
hyperbolic spaces existence of measures satisfying \eqref{eq:cac} (i.e.,
existence of invariant quasi-conformal streams in our terminology) was
established by Coornaert \cite{Coornaert93} (also by a generalization of
Patterson's construction).
\end{rem}

We are now ready to proceed to the main application of \thmref{th:equiv} which
is a description of the Hopf decomposition of (the measure class of) a
quasi-conformal stream invariant with respect to a discrete subgroup
$G\subset\Iso(\X)$. As the atomic part of such a stream is obviously
discontinual (so that by \thmref{th:main} its Hopf decomposition is completely
determined by the size of the point stabilizers), we can restrict our
considerations to the purely non-atomic case only.

\begin{thm} \label{th:main}
Let $G$ be a discrete group of isometries of a Gromov hyperbolic space $\X$.
Then the big horospheric limit set $\L_G^\hb$ is \textup{(mod~0)} the
conservative part of the boundary action of $G$ with respect to any purely
non-atomic $G$-invariant boundary quasi-conformal stream $\l$.
\end{thm}

\begin{proof}
Fix a reference point $o\in\X$ and consider the associated measure $\l_o$.
Without loss of generality we may assume that it is normalized, so that
$(\pt\X,\l_o)$ is a Lebesgue space. Any isometry of $\X$ has at most two fixed
boundary points, the group $G$ is countable, and the measure $m$ is purely
non-atomic. Therefore, the action of $G$ on the space $(\pt\X,\l_o)$ is free,
and we can apply \thmref{th:equiv}. By condition (iv) the orbit $G\o$ is
dissipative if and only if all the sets $\{g\in G: dg\l_o/d\l_o(\o)\ge t\}$
are finite, which, in view of \eqref{eq:cac}, is the same as the finiteness of
all the sets $\{g\in G: \b_\o(go,o)\ge t\}$, i.e., the same as the finiteness
of the intersection of the orbit $Go$ with any horoball centered at $\o$.
Thus, the orbit $G\o$ is dissipative if and only if it is contained in the
complement of $\L_G^\hb$.
\end{proof}

\begin{cor} \label{cor:main}
By \corref{cr:equiv}, the big horospheric limit set $\L_G^\hb$ ($\equiv$ the
conservative part of the boundary action) coincides with the divergence set of
the \emph{Poincar\'e--Busemann series}
$$
\left\{\o\in\pt\X: \sum_{g\in G} e^{D\b_\o(go,o)} = \infty \right\} \;.
$$
\end{cor}

\begin{rem}
For Fuchsian groups with respect to the visual stream on the boundary circle
\thmref{th:main} and \corref{cor:main} were proved by Pommerenke
\cite{Pommerenke76} (although in a somewhat different terminology, see the
discussion in \cite[Section 1]{Pommerenke82}). Pommerenke's argument uses
analytic properties of the Blaschke products and does not immediately carry
over to the higher dimensional situation. Sullivan \cite{Sullivan81} used a
more direct geometrical approach and established \thmref{th:main} for Kleinian
groups, again with respect to the visual stream (actually he considered the
small horospheric limit set $\L^\hs$ essentially showing that
$\L^\hb\setminus\L^\hs$ is a null set, see the Remark below). By extending
Sullivan's approach (with some technical complications) Tukia \cite{Tukia97}
proved \thmref{th:main} for Kleinian groups with respect to an arbitrary
invariant conformal stream.
\end{rem}

\begin{rem}
Our argument used the characterization of the conservative part of a free
action $G\act(X,m)$ as the set of all the points $x\in X$ such that for a
\emph{certain} $t>0$
$$
\{g\in G:dgm/dm(x)>t\} \quad \text{is infinite}
$$
(see condition (iv) from \thmref{th:equiv}), which in our setup is precisely
the big horospheric limit set $\L^\hb$. The small horospheric limit set
$\L^\hs$ corresponds to requiring that the above condition hold for \emph{any}
$t>0$, which, in general, is \emph{not} equivalent to conservativity (see
\remref{rem:special}). From this point of view the right object in the context
of studying conservativity of boundary actions is definitely $\L^\hb$ rather
than $\L^\hs$. The difference $\L^\hb\setminus\L^\hs$ is the set of all the
boundary points $\o\in\pt\X$ for which among the horoballs centered at $\o$
there are both ones containing \emph{finitely many points} from the orbit $Go$
of a fixed reference point $o\in\X$ and ones containing \emph{infinitely many
points} from $Go$. Sullivan \cite{Sullivan81} essentially proved that
$\L^\hb\setminus\L^\hs$ is a null set with respect to the visual stream for
Kleinian groups (also see the discussion of his result in \cite[Section
1]{Pommerenke82}); for subgroups of a free group (again with respect to the
uniform stream on the boundary of the ambient group) it was done in
\cite{Grigorchuk-Kaimanovich-Nagnibeda07}. We are not aware of any other
results of this kind; in particular, it is already not known for Kleinian
groups with respect to general invariant conformal streams, see
\cite{Tukia97}. Nonetheless, it seems plausible that $\L^\hb\setminus\L^\hs$
\emph{is a null set with respect to any invariant quasi-conformal stream on an
arbitrary Gromov hyperbolic space}.
\end{rem}

\bibliographystyle{amsalpha}
\bibliography{C:/Sorted/MyTex/mine}

\end{document}